\documentclass[12pt,leqno]{amsart}
\usepackage{amssymb}

\usepackage[all]{xy}
\usepackage{graphicx,amssymb,amsfonts,epsfig,amsthm,a4,amsmath,url}
\usepackage[latin1]{inputenc}

\swapnumbers

\newcounter{ictr}

\newcounter{actr}

\newtheorem{ThmIntro}{Theorem}
\newtheorem{CorIntro}[ThmIntro]{Corollary}
\newtheorem{PropIntro}[ThmIntro]{Proposition}

\newtheorem{thm}{Theorem}[section]
\newtheorem{cor}[thm]{Corollary}
\newtheorem{lem}[thm]{Lemma}
\newtheorem{clai}[thm]{Claim}
\newtheorem{prop}[thm]{Proposition}
\theoremstyle{definition}
\newtheorem{defn}[thm]{Definition}

\newtheorem{que}[thm]{Question}
\newtheorem{conj}[thm]{Conjecture}
\theoremstyle{remark}
\newtheorem{rem}[thm]{Remark}
\theoremstyle{cremark}

\newtheorem*{isop}{Isoperimetric profile}

\newtheorem*{asym}{Asymptotic behavior}
\newtheorem*{GES}{Geometrically Elementary Solvable groups}

\newtheorem*{isopb}{Isoperimetric profile inside balls}

\newtheorem{ex}[thm]{Example}

\numberwithin{equation}{section}

\newcommand{\Z}{\mathbf{Z}}

\newcommand{\N}{\mathbf{N}}
\newcommand{\F}{\mathbb{F}}
\newcommand{\R}{\mathbf{R}}
\newcommand{\C}{\mathbf{C}}
\newcommand{\Q}{\mathbf{Q}}

\newcommand{\BS}{\textnormal{BS}}
\newcommand{\Aff}{\textnormal{Aff}}

\newcommand{\bpr}{\noindent \textbf{Proof}: ~}

\newcommand{\epr}{~$\blacksquare$}

\newcommand{\Supp}{\textnormal{Supp}}

\title[Isoperimetric profile and random walks]{Isoperimetric profile and  random walks on locally compact solvable groups}
\author{Romain Tessera}

\date{\today}

\address{UMPA, ENS Lyon\\46 All\'ee d'Italie\\ 69364 Lyon Cedex\\ France}
\email{rtessera@umpa.ens-lyon.fr}

\begin{document}
\baselineskip=16pt

\maketitle

\begin{abstract}
We study a large class of amenable locally compact groups comprising all solvable algebraic groups over a local field and their discrete subgroups. We show that the isoperimetric profile of these groups is in some sense optimal among amenable groups. We use this fact to compute the probability of return of symmetric random walks, and to derive various other geometric properties.   
\end{abstract}


\section{Introduction}

We introduce a notion of large-scale foliation for metric measure
spaces and we prove that if $X$ is large-scale foliated by $Y$,
and if $Y$ satisfies a Sobolev inequality at large-scale, then so
does $X$. In particular the $L^p$-isoperimetric profile of $Y$
grows faster than the one of $X$. A special case is when $Y=H$ is
a closed unimodular subgroup of a locally compact unimodular group $X=G$. 
We apply this general observation to a special class of amenable groups, namely
geometrically elementary solvable groups (see the precise definition below). This class of amenable
locally compact groups is stable under quasi-isometries, and
contains all quotients of unimodular closed compactly generated
subgroups of the group of upper triangular matrices $T(d,k)$ for
every $d\in \N$ and every local field $k$. If $G$ is a
geometrically elementary solvable group with exponential growth,
we prove that the $L^p$-isoperimetric profile of $G$ satisfies
$j_{G,p}(t)\approx\log t$, for every $1\leq p\leq \infty$. As a
consequence, the probability of return of symmetric random walks
on such groups decreases like $e^{-n^{1/3}}.$ This work generalizes and unifies older results, where these estimates had been proved for unimodular amenable connected Lie groups  and their lattices \cite{P'}, and for some finitely generated groups such as the lamplighter or solvable Baumslag-Solitar's groups \cite{CGP}. 
Other  groups for which this behaviour has been established are finitely generated solvable torsion free groups with finite Prüfer rank \cite{PitSal}. 

We prove a stronger
statement when the group is a quotient of a solvable algrebraic group
over a $q$-adic field ($q$ is a prime), or a quotient of a closed subgroup of an amenable connected Lie group. Namely, for such  groups, the isoperimetric profile inside a ball 
grows linearly with the radius of the ball. This extends a previous result \cite{tessera}, where this had been shown for connected solvable  Lie groups (and their lattices).

The rest of the introduction is devoted to a more comprehensive and precise exposition of the results, together with some of their applications. In more details, after recalling a few definitions, we shall state our main results about the $L^p$-isoperimetric profile, and the $L^p$-isoperimetric profile in balls. Then we  shall explain how these estimates for the $L^2$-isoperimetric profile can be used to provide sharp lower bounds on the probability of return of random walks. Our estimates on the $L^p$-isoperimetric profile in balls will then be applied to compute the $L^p$-compression of these groups, and to prove that their first reduced $L^p$-cohomology vanishes for all $1<p<\infty$ (confirming a conjecture of Gromov). We then give a few interesting examples of groups in the class GES. We close the introduction with a short list of open problems.

\subsection*{Notation and basic definitions}

Let $G$ be a locally compact, compactly generated group equipped
with a left Haar measure $\mu$. Let $S$ be a compact
symmetric generating subset of $G$, i.e. $\bigcup_{n\in \N}S^n=G$.
Equip $G$ with the left-invariant word metric\footnote{To have a
real metric we must assume that $S$ is symmetric. However, this
assumption does not play any role in the sequel.} associated to
$S$, i.e. $d_S(g,h)=\inf\{n,\; g^{-1}h\in S^n\}$. The closed ball
of center $g$ and of radius $r$ is denoted by $B(g,r)$ and its
volume by $V(r)$.
Let $\lambda$ be the action of $G$ by left-translations on
functions on $G$, i.e. $\lambda(g)f(x)=f(g^{-1}x)$. Restricted to
elements of $L^p(G)$, $\lambda$ is called the left regular
representation of $G$ on $L^p(G)$.

For any $1\leq p\leq \infty$, and any subset $A$ of $G$, denote
$$J_p(A)=\sup_{f\in L^p(A)}\frac{\|f\|_p}{\sup_{s\in S}\|f-\lambda(s)f\|_p},$$
where $L^p(A)$ comprises all functions in $L^p(G)$ which are supported in $A$.

\begin{isop} (see for instance
\cite{C})
The $L^p$-isoperimetric profile, is the non-decreasing function
$$j_{G,p}(v)=\sup_{\mu(A)\leq v}J_p(A);$$
\end{isop} 
\begin{isopb}\cite{tessera}
The $L^p$-isoperimetric profile inside balls is the non-decreasing function
$$J^b_{G,p}(r)= J_p(B(1,r)).$$
\end{isopb}

\begin{asym}
Let $f,g: \R_+\to\R_+$ be two
monotonic functions. We write respectively $f\preceq g$, if there exists $C>0$ such that $f(t)=O(g(Ct))$ when $t\to \infty$. We write $f\approx g$ if both
$f\preceq g$ and $g\preceq f$. The asymptotic behavior of $f$ is
its class modulo the equivalence relation $\approx$.
\end{asym}

\subsection*{Main results}

In \cite{tes}, we proved that
\begin{ThmIntro}\label{QIThmIntro}
Let $(G,S)$ and $(H,T)$ be two compactly generated, locally
compact groups, equipped with symmetric generating subsets $S$ and
$T$ respectively. Then, the asymptotic behaviours of $j_{G,p}$,
$J^b_{G,p}$, for any $1\leq p\leq \infty$ does not depend on $S$.
Moreover, if $G$ and $H$ are both unimodular, and if $G$ is
quasi-isometric to $H$, then
$$j_{G,p}\approx j_{H,p},$$ and
$$J^b_{G,p}\approx J^b_{H,p}.$$
\end{ThmIntro}

Here we prove 
\begin{ThmIntro}\label{subgroupThmIntro}\textnormal{(see Corollary  \ref{groupisocor})}
Let $H<G$ be unimodular compactly generated groups and let
$1\leq p\leq \infty$. Then,
$$j_{G,p}\preceq j_{H,p}.$$
Moreover if $H$ is not distorted in $G$, then
$$J^b_{G,p}\preceq J^b_{H,p}.$$
\end{ThmIntro}
If the groups are finitely generated, these statements are much
easier to prove \cite{Ersch}.

We also show
\begin{PropIntro}\label{quotientThmIntro}\textnormal{(see Propositions \ref{quotientprop}, \ref{quotiencompacttprop} and \ref{cocompactprop})}
Let $1\to H\to G\to Q\to 1$ be an exact sequence of compactly generated locally compact groups.
Let $1\leq p\leq \infty$. Then, $j_{G,p}\preceq j_{Q,p},$ and
$J^b_{G,p}\preceq J^b_{Q,p}.$ Moreover if $H$ is compact, then  $j_{G,p}\approx j_{Q,p},$ and
$J^b_{G,p}\approx J^b_{Q,p}.$ If instead $Q$ is compact, then $j_{G,p}\approx j_{H,p},$ and
$J^b_{G,p}\approx J^b_{H,p}.$
\end{PropIntro}

According to a theorem of Coulhon and Saloff-Coste \cite{Coul}, if $G$ is a
compactly generated, locally compact group with exponential
growth, then $j_{G,p}(t)\preceq\log t$. On the other hand it is very
easy to see that $J_{G,p}(t)\leq 2t$. We shall now prove that
the converse inequalities hold for certain classes of
groups.  If $k$ is a field, let $T(d,k)$ be the group of invertible upper triangular matrices of size $d$ with coefficients in $k$.

\begin{GES}
The class GES of geometrically elementary solvable groups is the
smallest class of (compactly generated) locally compact groups
\begin{itemize}
\item comprising all unimodular closed compactly generated subgroups
of the group $T(d,k)$, for any integer $d$ and any local field
$k$; \item stable under taking finite products, quotients, and
unimodular closed compactly generated subgroups; \item stable
under quasi-isometry.
\end{itemize}
\end{GES}
As we will see below, this class contains as particular examples solvable Baumslag-Solitar groups, lamplighter groups, and polycyclic groups.  Let us now state our main result about this class.

\begin{ThmIntro}\label{localfieldThmIntro}
Let $G$ be a GES group. Then, for
every $1\leq p\leq \infty,$
$$j_{G,p}(t)\succeq\log t.$$
\end{ThmIntro}
This result was known for polycylic groups \cite{P},  connected amenable Lie groups \cite{P'}, 
for the lamplighter and other particular examples \cite{CGP}. To
prove Theorem~\ref{localfieldThmIntro}, we establish a stronger
result for the group of triangular matrices $T(d,k)$ over a local
field $k$, i.e. that $J^b_{G,p}(t)\succeq t$. 
Indeed, note that $J_{G,p}(r)\leq j_{G,p}(V(r)).$ So, in
particular, if the group has exponential growth,
$J^b_{G,p}(t)\succeq t$ implies $j_{G,p}\succeq \log t$.
The stability under
finite product being trivial, we obtain all geometric elementary
solvable groups using Theorems~\ref{QIThmIntro},
\ref{subgroupThmIntro} and \ref{quotientThmIntro}.

Restricting to groups with exponential growth, we obtain
\begin{CorIntro}
Let $G$ be an geometrically elementary solvable group with
exponential growth. Then, for every $1\leq p\leq \infty,$
$$j_{G,p}(t)\approx\log t.$$
\end{CorIntro}
It is natural to ask whether a geometrically elementary solvable group satisfies the stronger property 
$J^b_{G,p}(t)\succeq t$. Since its variable is the radius of a ball, $J^b_{G,p}$ is sensitive to distortion of the metric. This is reflected in the statement of Theorem \ref{subgroupThmIntro}, which shows that distortion is indeed the only limitation for such an improvement of our corollary. On the other hand, when dealing with $J^b_{G,p}$, the groups do not need to be unimodular (this assumption arises when we deal with volumes). 
\begin{ThmIntro}\label{padicThmIntro}
Let $k$ be a $q$-adic field. Let $G$ be either a quotient of a 
compactly generated algebraic group over $k$, or a quotient of a closed compactly generated subgroup of an almost connected  amenable Lie group.  Then, for every $1\leq p\leq \infty,$
$$J^b_{G,p}(t)\approx t.$$ In particular these groups have
controlled F\o lner sequences (see \cite{tessera'}).
\end{ThmIntro}
In \cite{tessera}, we proved it for connected amenable Lie groups,
lamplighter groups, and solvable Baumslag-Solitar groups. 

\subsection*{Application to random walks}

Let $G$ be a locally compact, compactly generated group. 
Using \cite[Theorem~9.2]{tes} (which
is a straightforward generalization of \cite[Theorem~7.1]{C}), we
obtain the following result.

\begin{ThmIntro}\label{randomwalkThmIntro}\textnormal{(see Theorem \ref{localfieldthm2})}
Let $G$ be a geometric elementary solvable group with exponential
growth, and let $\mu$ be a symmetric probability measure whose support is compact and generates $G$. Then, for every compact neighborhood $U$ of the neutral element, 
$$\mu^{2n}(U)\approx e^{-n^{1/3}}.$$
\end{ThmIntro}
As already mentioned, this fact was known for connected unimodular amenable Lie groups
\cite{P'}, for finitely generated torsion-free solvable groups
with finite Prüfer rank \cite{PitSal}, and for the lamplighter
group $F\wr \Z$, where $F$ is a finite group in \cite{CG}. 
Using a probabilistic approach, Mustapha \cite{Must'} was able to prove it 
for analytic $p$-adic unimodular groups (which are particular cases of geometric elementary solvable groups).

\subsection*{Application to $L^p$-compression}\label{compression}

\subsubsection*{Equivariant $L^p$-compression}

Recall that the equivariant $L^p$-compression rate of a locally
compact compactly generated group is the supremum of $0\leq \alpha
\leq 1$ such that there exists a proper isometric affine action
$\sigma$ on some $L^p$-spaces satisfying, for all $g\in G$,
$$\|\sigma(g).0\|_p\geq C^{-1}|g|_S^{\alpha}-C,$$ for some constant
$C<\infty$, $|g|_S$ being the word length of $g$ with respect to a
compact generating subset $S$.

It follows from \cite[Corollary~13]{tessera} that a group with
linear $L^p$-isoperimetric profile inside balls have equivariant
$L^p$-compression $B_p(G)=1$. Hence, we obtain
\begin{ThmIntro}
Let $k$ be a $q$-adic field. Let $H$ be  a closed
compactly generated subgroups of $T(d,k)$ whose Zariski closure
is compactly generated, and let $G$ be a quotient of  $H$. Then, $B_p(G)=1$ for any $1\leq p\leq
\infty$.
\end{ThmIntro}

\subsubsection*{Non-equivariant $L^p$-compression}

Recall that the $L^p$-compression rate of a metric space $(X,d)$
is the supremum of all $0\leq \alpha\leq 1$ such that there exists
a map $F$ from $X$ to some $L^p$-space satisfying, for all $x,y\in
X$,
$$C^{-1}d(x,y)^{\alpha}-C\leq \|F(x)-F(y)\|_p\leq d(x,y),$$
for some constant $C<\infty$.

Another theorem of Mustapha \cite{Must} says that an algebraic
compactly generated subgroup of $GL(d,k)$, where $k$ is a $q$-adic
field, is non-distorted in $GL(d,k)$. As $T(d,k)$ is co-compact in
$GL(d,k)$ and satisfies $B_p(T(d,k))=1$, we obtain
\begin{ThmIntro}
Let $k$ be a $q$-adic field. Let $G$ be an algebraic compactly
generated subgroups of $GL(d,k)$. Then, the $L^p$-compression rate
of $G$ satisfies $R_p(G)= 1$ for any $1\leq p\leq \infty$.
\end{ThmIntro}

\subsection*{Application to the first reduced $L^p$-cohomology}

Recall that a conjecture of Gromov \cite{G} states that all amenable discrete groups have trivial reduced $L^p$-cohomology, in degree at least $1$, and for any $1<p<\infty$. Cheeger and Gromov proved it when $p=2$ \cite{CG}. The question is obviously relevant in the more general setting of locally compact groups.  
In \cite{tessera'}, we were able to settle the conjecture in degree $1$ for the class of amenable groups with linear isoperimetric profile in balls. The present paper therefore extends the class of groups satisfying Gromov's conjecture in degree $1$. 

\begin{ThmIntro}
Let $G$ be as in Theorem \ref{padicThmIntro}. Then for every $1<p<\infty$,
$\overline{H}^1(G,\lambda_{G,p})=0$.
\end{ThmIntro}

\subsection*{Some remarks about the class GES}

As solvable connected Lie groups embed as a closed subgroups of $T(d,\C)$, we deduce from the discussion at the begining of Section \ref{GESsection} that unimodular amenable connected Lie groups (hence polycyclic groups) are in GES.

The class GES contains solvable Baumslag-Solitar groups $\BS(1,n)=\langle t,x\,|\;txt^{-1}=x^n\rangle,$ 
for any integer $n$ with $|n|\ge 1$. This group can be described as $\BS(1,n)=\Z[1/n]\rtimes\Z,$
where $\Z$ acts on $\Z[1/n]$ by multiplication by $n$. It has a matrix representation 
$$\BS(1,n)=\left\{\begin{pmatrix}
  n^k & P \\
  0 & 1 \\
\end{pmatrix}, k\in\Z, P\in \Z[1/n]\right\}.$$ Let $\Q_p$ denote the $p$-adic field, and define the ring $\Q_n$ as the direct product of all $\Q_p$ when $p$ ranges over the set of distinct prime divisors of $n$. Then the natural diagonal embedding of $\Z[1/n]$ into $\Q_n\oplus\R$ has discrete (and actually cocompact) image \cite[Chap.~IV, \S 2]{Weil}. Accordingly, $\BS(1,n)$ can be seen as a discrete subgroup of the product of two affine groups, namely $\Aff(\R)$ and $\Aff(\Q_n)$ (where $\Aff(k):=k\rtimes k^*$). 

Similarly if $p$ is a prime integer, the lamplighter group $\F_p\wr\Z=(\oplus_{\Z}\F_p)\rtimes \Z$ admits a matrix representation 
$$\F_p\wr\Z=\left\{\begin{pmatrix}
  X^k & P \\
  0 & 1 \\
\end{pmatrix}, k\in\Z, P\in \F_p[X,X^{-1}]\right\},$$ 
and therefore sits in $\Aff(\F_p((X)))\times\Aff(\F_p((X^{-1}))$ as a discrete subgroup. 

Note that the class GES also contains groups which are not virtually solvable. Indeed the lamplighter
group $F\wr \Z=F^{(\Z)}\rtimes \Z$, where $F$ is any finite group belongs to GES.
Namely, such a group is trivially quasi-isometric to any $F'\wr\Z$
where $F'$ has same cardinality as $F$  \cite{Di}. Hence one can take $F'$ to be a product of
$\F_q$ for finitely many primes $q$, which then implies that $F'\wr\Z$ is a subgroup
in a finite product of lamplighter groups $\F_q\wr \Z$.

Finally, let us mention that the class $GES$ contains finitely generated groups which are not residually finite (hence not linear) as shown by the following example due to Hall \cite{H}. Fix a prime
$q$ and consider the group of upper triangular $4$ by $4$
matrices:
$$G=\left\{\left(\begin{array}{cccc}
1 & 0 & x & z\\
0 & q^{n} & 0 & y \\
0 & 0 & q^{-n} & 0 \\
0 & 0& 0& 1
\end{array}\right); x,y,z\in Z[1/q]; n\in \Z \right\}.$$
Taking the quotient by the central infinite cyclic subgroup of
unipotent matrices $I+mE_{1,4}$ where $m\in \Z$, we obtain an
elementary solvable group which is non-residually finite since its
center is isomorphic to $\Z[1/q]/\Z$.

\subsection*{Open problems}

\begin{conj} We conjecture that all GES groups have controlled F\o lner pairs (see
Section~\ref{FolnerpairSection}), and therefore satisfy
$J_{G,p}(t)\approx t$. 
\end{conj}

\begin{que}
Is-it true that if $H<G$ are compactly generated locally compact amenable groups, $J_{G,p}(t)\leq J_{H,p}(t)$ (without the unimodularity assumption)? In particular, does $J_{G,p}(t)\approx t$
 hold for any closed subquotient of $T_d(A)$ with $A$ a finite product of local fields? 
\end{que}

\begin{que}
Are the classes GES, resp. of groups satisfying $J_{G,p}(t)\approx t$, stable under extension?
\end{que}

\begin{que} Does every group satisfying $j_{G,p}\succeq \log t$ belong to the class
GES? Or even better: is it always
quasi-isometric to a quotient of a closed subgroup of a linear group over a product of local fields? One can ask the same question for groups with equivariant Hilbert compression
rate $B(G)=1$.
\end{que}

\subsection*{Acknowledgments} I am grateful to Yves de
Cornulier for pointing me the group constructed by Hall \cite{H},
and for his useful remarks and corrections. 

\section{Organization of the paper}
\begin{itemize}
\item In Section~\ref{gradientDefSection}, we briefly recall the
notions of Sobolev inequalities at scale $h$ and the results of
\cite{tes} that we need here. We also discuss some subtleties happening when the group is not unimodular. 
\item In Section~\ref{foliationsection}, we prove our main technical result about
{\it large-scale foliations} of metric measure spaces. This is the main ingredient of the proof of Theorem \ref{subgroupThmIntro}.

\item In Section~\ref{HomogeneousSection}, we work out the case of
closed subgroups and quotients.

\item In Section~\ref{localfieldSection}, we prove
Theorems~\ref{localfieldThmIntro} and \ref{padicThmIntro}.

\item Finally, in Section~\ref{RandomSection}, we prove Theorem \ref{randomwalkThmIntro}. 
\end{itemize}

\section{Preliminaries: functional analysis at a given scale}\label{gradientDefSection}
The purpose of this section is to briefly recall the notions
introduced in \cite{tes}. By metric measure space $(X,d,\mu)$, we
mean a metric space $(X,d)$ equipped with a Borel
measure $\mu$  on $X$ such that bounded measurable subsets have finite measure. The volume of the closed ball
$B(x,r)$ is denoted by $V(x,r)$.

\subsection{The locally doubling property}

The metric measure spaces that we will consider satisfy a very
weak property of bounded geometry introduced in \cite{CS} in the
context of Riemannian manifolds.
\begin{defn}\label{doublinganyscaledefn}
We say\footnote{In \cite{CS} and in \cite{tess'}, the  local
doubling property is denoted $(DV)_{loc}$.} that a space $X$ is
locally doubling at scale $r>0$ if there exists a constant $C_r$
such that
$$\forall x\in X, \quad V(x,2r)\leq C_rV(x,r).$$
If it is locally doubling at every scale $r>0$, then we just say
that $X$ is locally doubling.
\end{defn}

\begin{ex}
Let $X$ be a connected graph with degree bounded by $d$, equipped
with the counting measure. The volume of balls of radius $r$
satisfies
$$\forall x\in X,\quad 1\leq V(x,r)\leq d^r.$$ In particular, $X$
is locally doubling.
\end{ex}

\begin{ex}
Let $(X,d,\mu)$ be a metric measure space and let $G$ be a locally
compact group acting by measure-preserving isometries. If $G$ acts
co-compactly, then $X$ is locally doubling.
\end{ex}

\subsection{Local norm of the gradient at scale
$h$}\label{largescalegradient}

The purpose of this section is to define a notion of ``local norm
of the gradient" (whose infinitesimal analogue is the modulus of
the gradient of a smooth function on a Riemannian manifold), which captures the local variations of a function defined on a metric measure space $(X,d,
\mu)$.  

For every $h>0$, we define an operator
$|\nabla|_{h}$ on $ L^{\infty}(X)$ by
$$\forall f\in  L^{\infty}(X),\quad |\nabla f|_{h}(x)=\sup\{|f(y)-f(x)|, d(x,y)\leq h\}.$$

\subsection{Sobolev inequalities}
Let $\varphi: \R_+\to\R_+$ be an increasing function and let $p\in
[1,\infty]$. The following formulation of Sobolev inequality was
first introduced in \cite{Cou} in the context of Riemannian
manifolds.

\begin{defn}\label{sobolevdef}
One says that $X$ satisfies a Sobolev inequality $(S_{\varphi}^p)$
at scale at least $h$ if there exist $C,C'>0$ such that
$$\| f\|_p\leq C\varphi(C'|\Omega|)\| |\nabla f|_{h}\|_p$$
where $\Omega$ ranges over all compact subsets of $X$, $|\Omega|$
denotes the measure $\mu(\Omega)$, and $f\in L^{\infty}(\Omega)$,
$ L^{\infty}(\Omega)$ being the set of elements of $L^{\infty}(X)$
with support in $\Omega$.
We say that $X$ satisfies a large-scale Sobolev inequality
$(S_{\varphi}^p)$ if it satisfies it at scale $h$ for $h$ large enough.
\end{defn}

\subsection{Isoperimetric profiles}\label{isoperimetricprofilesection}

Let $A$ be a measurable subset of $X$, and let $h>0$. We
define
$$J_{h,p}(A)=\sup_{f} \frac{\| f\|_p}{\||\nabla f|_{h} \|_p}$$
where the supremum is taken over functions $f\in L^{\infty}(A)$.

\begin{defn}\label{isoperimetricprofiledef}\cite{tes}
The $L^p$-isoperimetric profile $j_{X,h,p}$ (resp.
inside balls: $J^b_{X,h,p}$) is a nondecreasing function
defined by
$$j_{X,h,p}(v)=\sup_{|A|\leq v}J_{h,p}(A).$$
(resp. $J^b_{X,h,p}(t)=\sup_{x\in X}J_{h,p}(B(x,t)).$)

In the sequel, we will generally omit the scale, and only denote $j_{X,p}$ instead of $j_{X,h,p}$.
\end{defn}

\begin{rem} 
One can check that $j_{X,p}\preceq j_{X,q}$ is always
true when $p\leq q<\infty$. Moreover, in most cases (e.g. all
known examples of groups), $j_{X,p}\approx j_{X,q}$.
\end{rem} 

\subsection{Link between Sobolev inequalities and isoperimetric profiles}

Sobolev inequalities $(S_{\varphi}^p)$ can also be interpreted as
$L^p$-isoperimetric inequalities. Clearly, the space $X$ always
satisfies the Sobolev inequality $(S_{\varphi}^p)$ with $\varphi=
j_{X,p}$. Conversely, if $X$ satisfies $(S_{\varphi}^p)$ for a
function $\varphi$, then
$$j_{X,p}\succeq \varphi.$$

\subsection{$L^2$-profile and probability of return of random
walks}\label{probaSection}

The case $p=2$ is of particular interest as it contains some
probabilistic information on the space $X$.
Indeed, it was shown in  
\cite{CG} that for manifolds with bounded geometry, there is a good correspondence between 
upper bounds
of the large-time on-diagonal behavior of the heat kernel and  Sobolev inequality $(S_{\varphi}^2)$.
 In the survey \cite{C}, a similar statement is proved for the standard
random walk on a weighted graph. In \cite[Theorem~9.1]{tes}, we
give a discrete-time version of this theorem for general metric
measure spaces.

Let $(X,d,\mu)$ be a metric measure space. Consider a measurable family of probability measures $P=(P_x)_{x\in X}$ such that the operator on $L^2(X,\mu)$ defined by
$Pf(x)=\int_X f(y)dP_x(y)$
is self-adjoint. This is equivalent to say that the random walk with transition probabilities $P=(P_x)_{x\in X}$ is
{\it reversible} with respect to the measure $\mu$. We will also ask $P$ to focus on the geometry at scale $h>0$ in the following sense:  there exist a ``large" constant
$1\leq A<\infty$ and a ``small" constant $c>0$ such that for
($\mu$-almost) every $x\in X$:
\begin{itemize}
\item[(i)] $P_x\ll\mu;$

\item[(ii)] $p_x=dP_x/d\mu$ is supported in $B(x,Ah);$

\item[(iii)] $p_x$ is larger than $c$ on $B(x,h)$.
\end{itemize}
We will need the following particular case of
\cite[Theorem~9.2]{tes}.

\begin{thm}\label{random/sobolevThm}
Let $X=(X,d,\mu)$ be a metric measure space. Then, the large-scale
isoperimetric profile satisfies
$$j_{X,2}(t)\approx \log t,$$ if and only if for any reversible random
walk at scale large enough we have
$$\sup_{x\in X}p_x^{2n}(x)\approx e^{-n^{1/3}}\quad  \forall n\in \N.$$
\end{thm}

\subsection{Large-scale equivalence between metric measure
spaces}\label{largescaleDefSection}

\begin{defn}\label{largescaledef}
Let $(X,d,\mu)$ and $(X',d',\mu)$ two spaces satisfying the locally
doubling property. Let us say that $X$ and $X'$ are large-scale
equivalent if there is a function $F$ from $X$ to $X'$ with the
following properties
\begin{itemize}

\item[(a)] for every sequence of pairs $(x_n,y_n)\in
(X^2)^{\mathbb{N}}$
$$\left(d(F(x_n),F(y_n))\rightarrow\infty\right)\Leftrightarrow \left(d(x_n,y_n)\rightarrow\infty\right).$$

\item[(b)] $F$ is almost onto, i.e. there exists a constant $C$
such that $[F(X)]_C=X'$.

\item[(c)] For $r>0$ large enough, there is a constant $C_r>0$
such that for all $x\in X$
$$C_r^{-1}V(x,r)\leq V(F(x),r)\leq C_rV(x,r).$$
\end{itemize}
\end{defn}

\begin{rem}
Note that being large-scale equivalent is an equivalence relation
between metric measure spaces with locally doubling property.
\end{rem}

\begin{rem}
If $X$ and $X'$ are quasi-geodesic, then (a) and (b) imply that $F$
is roughly bi-Lipschitz: there exists $C\geq 1$ such that
$$C^{-1}d(x,y)-C\leq d(F(x),F(y))\leq Cd(x,y)+C.$$ 
In this case, (a) and (b) correspond to
the classical definition of a {\it quasi-isometry}.
\end{rem}

\begin{ex}\cite{tes} 
Consider the subclass of metric measure spaces including graphs with
bounded degree, equipped with the countable measure; Riemannian
manifolds with bounded geometry\footnote{Classically, this means: with Ricci curvature bounded from below, and with bounded radius of injectivity. These assumptions are generally used in order to perform some discretization of functional inequalities (see for instance \cite{Kan}).}, equipped with the Riemannian measure;
compactly generated, locally compact groups equipped with a left
Haar measure and a word metric associated to a compact, generating
subset. In this class, quasi-isometries are always large-scale
equivalences.
\end{ex}
 The following theorem generalizes \cite{CS}, which was established in the context of graphs and Riemannian manifolds with bounded geometry.

\begin{thm}\label{largescalethm}\cite[Theorem~8.1]{tes}
Let $F:X\rightarrow X'$ be a large-scale equivalence between two
spaces $X$ and $X'$ satisfying the locally doubling property.
Assume that for $h>0$ fixed, the space $X$ satisfies a Sobolev
inequality $(S^p_{\varphi})$ at scale $h$, then there exists $h'$,
only depending on $h$ and on the constants of $F$ such that $X'$
satisfies $(S_{\varphi}^p)$ at scale $h'$. In particular,
large-scale Sobolev inequalities are invariant under large scale
equivalence.
\end{thm}

\section{Large-scale foliation of a metric measure
space and monotonicity of the isoperimetric
profile}\label{foliationsection}

\begin{defn}\label{foliationdef}
Let $X=(X,d_X,\mu)$ and $Y=(Y,d_Y,\lambda)$ be two metric measure
spaces satisfying the locally doubling property. We say that $X$ is
large-scale foliated (resp. normally large-scale foliated) by $Y$ if
it admits a measurable partition $X=\sqcup_{z\in Z}Y_z$ satisfying
the following two first (resp. three)  conditions.
\begin{itemize}
\item(measure decomposition) There exists a measure $\nu$ on $Z$ and a measure
$\lambda_z$ on $\nu$-almost every $Y_z$ such that for every
continuous compactly supported function $f$ on $X$,
$$\int_X f(x)d\mu(x)=\int_Z\left(\int_{Y_z}f(t)d\lambda_z(t)\right)d\nu(z).$$
The subsets $Y_z$ are called the leaves, and the space $Z$ is
called the base of the foliation. 

\item(large-scale control) For $\nu$-almost
every $z$ in $Z$, there exists $\alpha_z>0$ and a large scale equivalence $h_z: (Y,d_Y,\lambda) \to (Y_z,d_X,\alpha_z\nu_z)$, which is uniform with respect to $z\in Z$. 

\item(measure normalization)
There exists a constant $1\leq C<\infty$ such that
for every $z\in Z$ and every $x\in Y_z$,
$$C^{-1}\leq V_{X}(x,1),  V_{Y_z}(x,1)\leq C.$$
In particular, we can take $\alpha_z=1$.
\end{itemize}
\end{defn}

Recall that the compression of a map $F$ between two metric space
$X$ and $Y$ is the function $\rho$ defined by
$$\forall t>0,\quad \rho(t)=\inf_{d_X(x,x')\geq t}d_Y(F(x),F(x')).$$
\begin{defn}
We call the compression of a large-scale foliation of $X$ by $Y$ the
function $$\rho(t)=\inf_{z\in Z}\rho_z(t)$$ where $\rho_z$ is the
compression function of the large-scale equivalence $h_z$.
\end{defn}
A crucial example that we will consider in some details in the next
Section~ is the case when $Y=H$ is a closed subgroup of a locally
compact group $G=X$ such that $G/H$ carries a $G$-invariant measure.

\begin{thm}\label{foliationprop}
Let $X=(X,d,\mu)$ and $Y=(Y,\delta,\lambda)$ be two metric measure
spaces satisfying the locally doubling property. Assume that $X$
is normally large-scale foliated by $Y$. Then if $Y$ satisfies a
Sobolev inequality $(S_{\varphi}^p)$ at scale $h$, then $X$
satisfies $(S_{\varphi}^p)$ at scale $h'$, for $h'$ large enough.
In other words, if $j_{X,p}$ and $j_{Y,p}$ denote respectively the
$L^p$-isoperimetric profiles of $X$ and $Y$ at scale $h$ and $h'$,
then $$j_{Y,p}\succeq j_{X,p}.$$ Moreover, if $\rho$ is the
compression of the large-scale equivalence, then
$$J_{Y,p}^b\succeq J_{X,p}^b\circ \rho.$$
The latter result is true under the weaker assumption that $X$ is
merely large-scale foliated by $Y$.
\end{thm}

The main difficulty comes from the fact that we need to control
the measure of the support of the restriction to a leaf of a
function defined on $X$. On the contrary,  the control on the diameter of the support follows trivially from   the definition of
$\rho$ (and does not require the normalization condition). Hence the second inequality
($J_{Y,p}^b\succeq J_{X,p}^b\circ \rho$) is much easier than the first one, and therefore left to the
reader.

\begin{defn}
A subset $A$ of a metric space is called $h$-thick if it is a
reunion of closed balls of radius $h$.
\end{defn}

The following lemma implies that we can restrict
 to functions with thick support.

\begin{lem}\label{thickprop}(see \cite[Proposition~8.3]{tes})
Let $X=(X,d,\mu)$ be a metric measure space. Fix some $h>0$ and
some $p\in [1,\infty]$. There exists a constant $C>0$ such that
for any $f\in   L^{\infty}(X)$, there is a function $\tilde{f}\in
L^{\infty}(X)$ whose support is included in a $h/2$-thick subset
$\Omega$ such that
$$\mu(\Omega)\leq \mu(\Supp(f))+C$$
and for every $p\in [1,\infty]$,
$$\frac{\| |\nabla \tilde{f}|_{h/2} \|_p}{\|\tilde{f}\|_p}\leq C\frac{\| |\nabla f|_{h}\|_p}{\| f\|_p}.$$
\end{lem}

On the other hand, the locally doubling property ``extends" to
thick subsets in the following sense.

\begin{lem}\label{DVthickprop}
Let $X$ be a metric measure space satisfying the locally doubling
property. Fix two positive numbers $u$ and $v$. There exists a
constant $C=C(u,v)<\infty$ such that for any $u$-thick subset
$A\subset X$, we have
$$\mu([A]_v)\leq C\mu(A),$$
where $[A]_v:=\{x\in X,\; d(x,A)\leq v\}.$
\end{lem}
\bpr
Since $A$ is $u$-thick, it can be written as a union $\bigcup_{i\in I}B(x_i,u)$. Take a maximal subset $J\subset I$ such that the balls $B(x_j,u)$, with $j\in J$ are disjoint. By maximality of $J$, $A\subset  [A]_v\subset \bigcup_j B(x_j,4u+v)$. Indeed, let $x\in [A]_v$, then the ball $B(x,v +3u)$ contains a ball $B(x_i,u)$ for some $i\in I$. Therefore maximality implies that this ball must meet $B(x_j,u)$ for some $j\in J$, which implies the statement.  The lemma now follows from the locally doubling condition.
\epr

\

Finally, we will  need
\begin{lem}\label{leafLemma}
Assume that $X$ is normally large-scale foliated by $Y$. For every
$z\in Z$, let $[Y_z]_1$ be the $1$-neighborhood of $Y_z$ in $X$.
The inclusion map $(Y_z,d_X, \lambda_z) \to ([Y_z]_1,d_X,\mu)$ is a large-scale equivalence,
uniformly w.r.t. $z$.
\end{lem}
\noindent{\bf Proof of Lemma~\ref{leafLemma}.} The two metric
conditions (a) and (b) for being a large-scale equivalence (see
Definition~\ref{largescaledef}) are trivially satisfied here. It
remains to compare the volume of balls of fixed radius. But this
follows from the third condition of Definition~\ref{foliationdef} and the locally doubling property of the spaces $X$ and $Y$.\epr

\

\noindent{\bf Proof of Theorem~\ref{foliationprop}.} Throughout the
proof, we will use the letter $C$ as a generic constant, which might possibly take different values. Assume that
$Y$ satisfies the Sobolev inequality $(S_{\varphi}^p)$. Let
$\Omega$ be a compact subset of $X$ and $f\in L^{\infty}(\Omega)$.
We want to prove that there exists $h'$, depending only on $h$ and the spaces such that $f$ satisfies $(S_{\varphi}^p)$ at
scale $h'$. By Lemma~\ref{thickprop}, we can assume that $\Omega$
is $1$-thick. For every $z\in Z$, denote by $f_z$ the restriction
of $f$ to $Y_z$ and $\Omega_z=\Omega\cap Y_z$.

\begin{clai}
There exists $C<\infty$ such that for every $z\in Z$
$\lambda_z(\Omega_z)\leq C\mu(\Omega)$.
\end{clai}
\bpr  Lemma \ref{leafLemma} and \cite[Proposition 8.5]{tes} imply that 
$ \lambda_z(\Omega_z)$ is less than a constant times $\mu([\Omega_z]_1)$, which is obviously less than $\mu([\Omega]_1)$, which by Lemma \ref{DVthickprop} is less than a constant times $\mu(\Omega)$. 
 \epr

\

By Theorem~\ref{largescalethm}, there exists $h'>0$ such that
$Y_{z}$ satisfies $(S_{\varphi}^p)$ at scale $h'$, uniformly with
respect to $z\in Z$. So for every $z\in Z$,
$$\| f_z\|_p\leq C\varphi(C\lambda_z(\Omega_z))\||\nabla f_z|_{h'}\|_p.$$
Since $\lambda_z(\Omega_z)\leq C\mu(\Omega)$ and $\varphi$ is
nondecreasing, we have
$$\| f_z\|_p\leq C\varphi(C\mu(\Omega))\||\nabla f_z|_{h'}\|_p.$$
Moreover, we have
$$\| f\|_p^p =\int_{Z}\| f_z\|_p^pd\nu(z).$$
Clearly, since $Y_z$ is equipped with the induced distance, for
every $z\in Z$ and every $x\in Y_z$,
$$|\nabla f|_{h'}(x)\geq |\nabla f_z|_{h'}(x).$$ Therefore,
$$\| |\nabla f|_{h'}\|_p^p\geq \int_Z\||\nabla f_{z}|_{h'}\|_{p}^pd\nu(z).$$
We then have
$$\| f\|_p\leq C\varphi(C\mu(\Omega))\||\nabla f|_{h'}\|_p,$$
and we are done. \epr

\section{Application to locally compact groups}\label{HomogeneousSection}

\subsection{The case of groups: left and right translations}\label{unimodularity}

Let $G$ be a locally compact, compactly generated group, and let
$S$ be a generating set. Let $g\in G$ and let $f\in L^p(G)$ for
some $1\leq p\leq \infty.$ We have
$$|\nabla f|_{1}(g)=\sup_{s\in S}|f(gs)-f(g)|.$$
In other word, if $\rho$ is the action of $G$ by right-translation
on functions, i.e. $\rho(g)f(x)=f(xg)$, the isoperimetric profiles
are therefore given by
$$j_{G,p}(m)=\sup_{|\Supp(f)|\leq m}\frac{\|f\|_p}{\|\sup_{s\in S}|f-\rho(s)f|\|_p},$$
and 
$$J^b_{G,p}(r)=\sup_{\Supp(f)\subset B(1,r)}\frac{\|f\|_p}{\|\sup_{s\in S}|f-\rho(s)f|\|_p}.$$

\begin{prop}
We have $$j_{G,p}\approx\sup_{|\Supp(f)|\leq m}\frac{\|f\|_p}{\sup_{s\in S}\|f-\rho(s)f\|_p},$$ and
$$J^b_{G,p}(r)\approx\sup_{\Supp(f)\subset B(1,r)}\frac{\|f\|_p}{\sup_{s\in S}\|f-\rho(s)f\|_p}.$$
\end{prop}
\bpr Let us prove it for $j_{G,p}$ (the proof for $J^b_{G,p}$ is similar)

Note that $$\left\|\frac{1}{|S|}\left(\int_{S}|f-\rho(s)f|^pds\right)^{1/p}\right\|_p^p= \frac{1}{|S|}\int_{S}\|f-\rho(s)f|\|_p^pds\leq \sup_{s\in S}\|f-\rho(s)f\|_p^p\leq \|\sup_{s\in S}|f-\rho(s)f|\|_p^p,$$ we get

\begin{equation}\label{eq:leq}
j_{G,p}\leq \sup_{|\Supp(f)|\leq m}\frac{\|f\|_p}{\sup_{s\in S}\|f-\rho(s)f\|_p}\leq \sup_{|\Supp(f)|\leq m}\frac{\|f\|_p}{\left\|\frac{1}{|S|}\left(\int_{S}|f-\rho(s)f|^pds\right)^{1/p}\right\|_p}.
\end{equation} 
But in \cite[Section 7]{tes}, it was proved that 
\begin{equation}\label{eq:approx}
j_{G,p}\approx \sup_{|\Supp(f)|\leq m}\frac{\|f\|_p}{\left\|\frac{1}{|S|}\left(\int_{S}|f-\rho(s)f|^pds\right)^{1/p}\right\|_p}
\end{equation}
So combining (\ref{eq:leq}) and (\ref{eq:approx}) proves the proposition. \epr

\subsection*{Left or right?}
One may wonder why we chose to define the isoperimetric
profiles with left-translations in the introduction, since according to the previous proposition, it seems that the correct definition should be with right-translations. Indeed, since the metric is left-invariant, we have for all $s\in S$, $d(g,gs)=1$ whereas in general $d(s^{-1}g,g)$ is not even bounded. Here are the reasons of this
choice
\begin{itemize}
\item[-] if the group $G$ is unimodular, then the isoperimetric
profiles {\it are the same}, whether we define them by left-translations, or by
right-translations;

\item[-] Suppose that the group $G$ is non-unimodular. Then if we define it with 
left-translations and if the group is amenable, 
$j_{G,p}=\infty$. On the other hand if we define them with
right-translations, then both $j_{G,p}$ and $J^b_{G,p}$ are bounded and therefore behave as for a non-amenable group.

\item[-]  Finally,  if we define $J^b_{G,p}$ with
left-translations and if the group is amenable, then it is always a non-bounded increasing function.
\end{itemize}
In conclusion, we see that if the group is non-unimodular, then the asymptotic behavior of $j_{G,p}$ does not contain any interesting information on the group, whereas $J^b_{G,p}$ is interesting only if it is defines with left-translations. 

In the following sections, we will not change our notation but rather indicate
whether we consider a ``left-profile" or a ``right-profile" on
$G$.

\subsection{Closed subgroups}

\begin{prop}\label{subgroupfoliation}
Let $H$ be a closed compactly generated subgroup of a locally
compact compactly generated group $G$. Assume that the quotient
$G/H$ carries a $G$-invariant Borel measure, then $G$ is normally
large-scale foliated by $H$.
\end{prop}

\bpr  Let $\nu$ be a $G$-invariant $\sigma$-finite measure on the
quotient $Z=G/H$. Since $\nu$ is $G$-invariant, up to normalize it,
one can assume that for every continuous compactly supported
function $f$ on $G$,
$$\int_G f(g)d\mu(x)=\int_Z \left(\int_{H}f(gh)d\lambda(h)\right)d\nu(gH).$$
We claim that the partition $G=\sqcup_{gH\in Z} gH$ satisfies the third conditions of  Definition~\ref{largescaledef}.
Clearly, the first one follows from the above decomposition of $\mu$. For
every $g\in G$, the left-translation by $g$ is an isometry on $G$, so that all the left cosets of $H$ are isometric (we consider these cosets equipped with the distance induced from the one of $G$).
On the other hand, since $H$ is a closed subgroup, the inclusion
map $H\to G$ is a uniform embedding, i.e. satisfies (a)
of Definition~\ref{largescaledef}. This proves the second condition. Finally, the last
condition follows from the left-invariance of both $\nu$ and
$\mu$. Namely, the left-invariance of $\nu$ implies that, for
every $g\in G$, the mesure $\lambda_g$ on $gH$ is the image
of $\lambda$ under the map $h\to gh$ from $H$ to $gH$. \epr

In \cite[Lemma 4]{Ersch}, it is proved that if $H$ is finitely
generated subgroup of a finitely generated group $G$, then
$j_{H}\preceq j_{G}$. Here is a generalization of this easy result.
\begin{cor}\label{groupisocor}
Let $H$ be a closed, compactly generated subgroup of $G$ and let
$1\leq p\leq \infty$. Assume that $G/H$ carries a $G$-invariant
measure. Then,
\begin{itemize}
\item The right-profiles satisfy $j_{G,p}\preceq j_{H,p};$ \item A
weaker conclusion holds for the right-profiles inside balls:
$J^b_{G,p}\preceq J^b_{H,p}\circ\rho,$ where $\rho$ is the
compression of the injection $H\hookrightarrow G.$
\end{itemize}
\end{cor}
\bpr
This follows from Theorem \ref{foliationprop}  and Proposition \ref{subgroupfoliation}.
\epr

\begin{rem}\label{remunimodular}
Corollary~\ref{groupisocor} holds in particular when $G$ and $H$ are
both unimodular. Actually this is the only interesting situation
since, by \cite[Lemma~11.10]{tes}, a non-unimodular group 
always has a bounded right-profile, i.e.  $\||\nabla f|_h\|_p\geq c\|f\|_p$ for some $c>0$ depending only on 
$p\geq 1$ and $h\geq 1$. On the other hand, if $H$ is
non-unimodular and if $G$ is unimodular and amenable, then, by
\cite[Proposition~11.11]{tes} all the conclusions of
Corollary~\ref{groupisocor} fall appart\footnote{For example,
consider the non-unimodular group $H$ of positive affine
transformations of $\R$: this group, equipped with its
left-invariant Riemannian metric is isometric to the Hyperbolic
plane. In particular, it has a bounded isoperimetric right-profile. On
the other hand, it is a closed subgroup of the solvable unimodular
Lie group Sol, whose isoperimetric profile $j_{G,p}$ is
asymptotically equivalent to $\log t$ \cite{P'}.}.
\end{rem}

\subsection{Quotients and co-compact subgroups}

\begin{prop}\label{quotientprop}
Let $Q=G/H$ be the quotient of a locally compact, compactly
generated group $G$ by a closed normal subgroup $H$. Then for all
$1\leq p\leq \infty$, the left-profiles satisfy $j_{G,p}\preceq
j_{Q,p}$ and the left-profiles in balls satisfy $J^b_{G,p}\preceq
J^b_{Q,p}.$
\end{prop}
\bpr We denote by $\pi$ the projection on $G/H$. Let us equip $G$
and $H$ with left Haar measures $\mu$ and $\nu$. Take a Haar measure
$\lambda$ on $H$ such that for every continuous compactly
supported function $f$ on $G$,
$$\int_G f(g)d\mu(g)=\int_Q \left(\int_{H}f(gh)d\lambda(h)\right)d\nu(gH).$$
Let $S$ be a symmetric compact generating subset of $G$ and
consider its image $T$ by the projection on $Q$. The projection
$\pi$ is therefore $1$-Lipschitz between $(G,S)$ and $(Q,T)$. For
every $1\leq p<\infty$, consider the application $\Psi: \;
C_0(G)\to C_0(Q)$ defined by
$$\Psi(f)(gH)=\left(\int_H |f(gh)|^pd\lambda(h)\right)^{1/p}.$$
Clearly, the support of $\Psi(f)$ is the projection of the support
of $f$. Moreover, $\Psi$ preserves the $L^p$-norm. Take $s\in S$
and $t=\pi(s)$, $g\in G$ and $q=\pi(g)$,
\begin{eqnarray*}
|\Psi(f)(t^{-1}q)-\Psi(f)(q)| & = & \left(\int_H
|f(s^{-1}gh)|^pd\lambda(h)\right)^{1/p}-\left(\int_H
|f(gh)|^pd\lambda(h)\right)^{1/p}\\
                              & \leq & \left(\int_H
|f(s^{-1}gh)-f(gh)|^pd\lambda(h)\right)^{1/p}.
\end{eqnarray*}
Therefore,
$$\|\lambda(t)\Psi(f)-\Psi(f)\|_p\leq \|\lambda(s)f-f\|_p,$$
and we are done. \epr

\begin{prop}\label{quotiencompacttprop}
Let $Q=G/K$ be the quotient of a locally compact, compactly
generated group $G$ by a closed normal compact subgroup $K$. Then for all
$1\leq p\leq \infty$, the (left or right) profiles satisfy $j_{G,p}\approx
j_{Q,p}$ and the (left or right) profiles in balls satisfy $J^b_{G,p}\approx
J^b_{Q,p}.$
\end{prop}
\bpr
Note that by the previous proposition, we only need to show that $j_{G,p}\succeq
j_{Q,p}$ and  $J^b_{G,p}\succeq
J^b_{G,p}.$ Consider compact generating sets $S$ and $\pi(S)$ on respectively $G$ and $Q$, and suppose that the Haar measures on $G$ and $Q$ are such that $\mu(\pi^{-1}(A))=\nu(A)$ for every Borel  subset $A\subset Q$ with finite measure. Note that the composition with $\pi$ defines a map $C_0(Q)\to C_0(G)$ which induces an isometry $\Phi: L^p(Q)\to L^p(G)$. On the other hand since $K$ is normal, one checks easily that $\Phi$ commutes with both the left and right regular representations. Indeed, 
$\lambda(s)\Phi(f)=\Phi(\lambda(\pi(s))f)$ (and same for $\rho$). Combining this with the fact that $\Phi$ is an isometry, yields the following equality. 
For all $s\in S$, 
$$\|\lambda(s)\Phi(f)-\Phi(f)\|_p= \| \lambda(\pi(s))f-f\|_p.$$
This proves the proposition.
\epr

\begin{prop}\label{cocompactprop}
Let $1\to H\to G\to Q\to 1$ be a short exact sequence of compactly generated locally compact groups such that $H$ is closed and $Q$ is compact. Then for all
$1\leq p\leq \infty$, the profiles satisfy $j_{G,p}\approx
j_{H,p}$ and the  profiles in balls satisfy $J^b_{G,p}\approx
J^b_{H,p}.$
\end{prop}
\bpr 
We only need to show the inequalities $j_{G,p}\succeq
j_{H,p}$ and  $J^b_{G,p}\succeq J^b_{G,p}.$  Let $\eta: Q\to G$ be a measurable section such that $K=\eta(Q)$ is relatively compact. Each element $g\in G$ can be uniquely  written as $\eta(q)h$, with $q\in Q$ and $h\in H$, and the left Haar measure on $G$ identifies (up to normalization) as the direct product\footnote{To see this, observe that if $g\in G$, $q\in Q$ and $h\in H$, then $g\eta(q)h=\eta(q')th$, where $q'=\bar{g}q\in Q$, and $t=\eta(q')^{-1}s\eta(q)\in H$.} of the left Haar measures of $Q$ and $H$.   Define an isometric embedding $\Phi: L^p(H)\to L^p(G)$ by 
$$\Phi(f)(kh)=f(h).$$ Suppose that the Haar measure of $Q$ has total mass equal to 1. Then the measure of the support of $\Phi(f)$ equals the measure of the support of $f$. On the other hand, since the inclusion $H\to G$ is a quasi-isometry, the diameter of the support of $\phi(f)$ is less than a constant times the diameter of the support of $f$. 

Finally let $S$ be a compact generating subset of $G$. For all $s\in S$, $h\in H$, $q\in Q$ and $f\in L^p(H)$, we have
$$\Phi(f)(s\eta(q)h)-\Phi(f)(\eta(q)h)=\Phi(f)(\eta(q')th)-\Phi(f)(\eta(q)h)=
f(th)-f(h),$$
where $q'=\bar{g}q\in Q$, and $t=\eta(q')^{-1}s\eta(q)$ is an element of  $K^{-1}SK\cap H$ which is relatively compact in $H$. 
The proposition now follows trivially. 
\epr

\section{Geometrically elementary solvable groups}\label{localfieldSection}\label{GESsection}
As already mentioned in the introduction, Theorems~\ref{localfieldThmIntro} and
\ref{padicThmIntro} will be consequences of the following facts: 

\begin{itemize}
\item the stability
results: Corollary~\ref{groupisocor}, Propositions~\ref{quotientprop}, \ref{quotiencompacttprop} and \ref{cocompactprop}. 
\item the result of Mustapha \cite{Must} (for Theorem \ref{padicThmIntro}),
recalled in the introduction, 
\item the case
of algebraic groups over p-adic fields (and $T(d,k)$ for any local field), treated in the last subsection.
\end{itemize}

Let us now briefly explain how to deduce Theorem \ref{padicThmIntro}. 
Recall that every connected amenable Lie group has a closed, connected, co-compact, normal solvable subgroup (namely its solvable radical).
Moreover by a result of Mostow \cite{M},  any closed compactly generated subgroup $G$ of a solvable connected Lie group has a normal compact subgroup $K$, such that $G/K$ embeds as  co-compact (hence non-distorted) subgroup of a solvable connected Lie group.

\subsection{F\o lner pairs and isoperimetric profile}\label{FolnerpairSection}
The notion of F\o lner pairs was introduced in \cite{BCG} in order to produce lower bounds on the probability of return of symmetric random walks on graphs (or of the Heat kernel on a Riemannian manifold). In \cite{tessera}, we defined a slightly different notion, called controlled F\o lner pairs.  
\begin{defn}
Let $G$ be a locally compact compactly generated group and let $S$
be a compact symmetric generating subset of $G$. A sequence
$(F_n,F_n')$ of pairs of compact subsets is a sequence of
controlled F\o lner pairs if there is a constant $C<\infty$ such
that for all $n\in \N$,
\begin{itemize}
\item $\mu(F_n')\leq C\mu(F_n);$ \item $S^nF_n\subset F_n'$.
\end{itemize}
\end{defn}

We will need the following easy fact.

\begin{prop}\label{controlfolnerProp}\cite[Proposition~4.9]{tessera}
If $G$ has a sequence of controlled F\o lner pairs, then for all
$1\leq p\leq \infty,$ $J^b_{G,p}(t)\succeq t.$ Moreover, if $G$ is
unimodular, then
$j_{G,p}(t)\succeq \log t.$
\end{prop}

\subsection{Algebraic groups over a p-adic field}
Recall that an algebraic group over $k$ is isomorphic to a semi-direct product $T\ltimes U$ where $T\simeq  (k^*)^d$ acts semi-simply on the unipotent radical $U$. Moreover, $G$ is compactly generated if and only if this action is non-degenerate in the sense that the weights of the action of $T$ on $G/[G,G]$ are non-zero (see \cite[Theorem 13.4]{BT} or \cite[Theorem 3.2.3]{Ab}).
\begin{thm}
Let $k$ be a non-archimedean local field, and let $G= T\ltimes U$ be a compactly generated  algebraic group over $k$, with $T\simeq  (k^*)^d$ acting semi-simply on the unipotent radical $U$. Then $G$ has a sequence of F\o
lner pairs. In particular (see
Proposition~\ref{controlfolnerProp}), it satisfies
$J^b_{G,p}(t)\approx t$.
\end{thm}

\bpr For the sake of concreteness, we will write the proof in the special case where $G=T(d,k)$ is the group of invertible upper triangular matrices of size $d$ with coefficients in $k$.

Let us assume that $d\geq 2$ (the case $d=1$ being trivial). Let
$v$ be a discrete valuation on $k$ (with values in $\Z$), and for every $x\in k$, let
$|x|=e^{-v(x)}$ be the corresponding norm. We have
$$|x+y|\leq \max\{|x|,|y|\},$$
and
$$|xy|=|x||y|.$$
Let $k_n$ be the (compact) subring of $k$ consisting of elements
$y\in k$ of norm $|y|\leq e^n$. We fix a uniformizer $x_0\in k$ (i.e. such that
$v(x_0)=-1$). We have
$$x_0k_n=k_{n+1}.$$

\medskip

Let $U$ be the subgroup of $G$ consisting of unipotent elements,
and let $T\simeq (k^*)^d$ be the subgroup of diagonal elements. We
have a semidirect product
$$G=T\ltimes U.$$
For every $n\in \N$, let $U_n$ be the compact normal subgroup of
$U$ consisting of unipotent matrices such that for $1\leq i<j\leq d$, the $(i,j)$'th coefficient lies in $k_{(j-i)n}$. 
We also consider the compact subset $T_n$ of $T$ defined by
diagonal matrices whose diagonal coefficients and their inverses have norm less or equal than $e^n$.

Let us identify $G$ with the cartesian product $T\times U$, where
the group law is given by
$$(t,u)(s,v)=(ts,u^sv),$$
where $u^s=s^{-1}us.$ We define a compact subset $S$ of $G$ by
$$S=T_1\cup U_0.$$
Let $t\in T_1.$ An easy computation shows
that for every $n\in \N,$
\begin{equation}\label{contractibilityEq1}
 t^{-1}U_n t\subset U_{n+1}.
\end{equation}
Moreover, $t_0=(x_0^{d-1},x_0^{d-2},\ldots,1)\in T_{d-1}$ satisfies, for every $n\in \N,$
\begin{equation}\label{contractibilityEq2}
 t_0^{-1}U_n t_0\subset U_{n+1}.
\end{equation}

Note that since $G=\bigcup_n T_n\times U_n$, this implies that $S$ is a generating subset of
$G$. 
On the other hand, we deduce from (\ref{contractibilityEq2}) that 
$$U_{n} \subset S^{2dn+1}.$$
As $T_n\subset S^n$, we have
$$T_n\times U_n\subset S^{2dn+n+1}.$$
\begin{clai} For all $n\geq 1$, $S^n\subset T_n\times
U_{n}.$

\end{clai} 
\bpr This is true for $n=1$. Now, assume that this is true for $n\geq 1$, and
take an element $g=(t,u)$ in $T_{n}\times U_{n},$ and an
element $h$ of $S$. Let us check that $gh\in T_{n+1}\times
U_{n+1}$. First, assume that $h=(s,1)\in T_1$. Then,
$$gh=(ts,u^s)\in T_{n+1}\times s^{-1}U_{n}s^{-1}$$
By (\ref{contractibilityEq1}), $ T_{n+1}\times s^{-1}U_{n}s^{-1}\subset T_{n+1}\times
U_{n+1}.$ 

Now, if $h=(1,v)\in U_0$, then 
$gh=(t,uv)\subset T_n\times U_{n}.$ \epr

\

Now, let $F_n=T_n\times U_{2n}$ and $F'_n=T_{2n}\times U_{2n}$. We
claim that $(F_n,F'_n)$ is a sequence of F\o lner pairs. As
$F'_n\subset S^{4dn+2n+1}$ and $|F'_n|=2|F_n|$, we just need to check
that $S^nF_n\subset F'_n.$ Let $g=(t,u)\in S^n\subset T_n\times
U_{n}$ and $g'=(s,v)\in F_n$. By an immediate induction, (\ref{contractibilityEq1}) implies that
$s^{-1}U_ns^n\subset U_{2n}.$ Hence,
$$(t,u)(s,v)=(ts,u^sv)\in T_{2n}\times U_{2n}U_{n}=T_{2n}\times U_{2n},$$
which finishes the proof.\epr

\section{Random walks}\label{RandomSection}

We will prove a slightly more general result than the one stated in the introduction.
Let $G$ be a locally compact, compactly generated group and let $S$ be a compact generating set of $G.$ Let
$(X,\mu,d)$ be  a measure space, equipped with a measurable metric, on which $G$ acts measurably by isometries, and such that for every $x\in X$ the orbit map from $G$ to $X$: $g\to gx$ is a large-scale equivalence. If $X$ satisfies these asumptions, we will call it a geometric $G$-space. A typical example is if $G$ acts continuously, properly and co-compactly on $X$.

Let $x_0$ be a point in $X$, and $r_0$ be such that every $x\in X$ lies at distance $<r_0$ from the orbit $Gx_0$. 

For every $x\in X$, let $\nu_x$ be a probability measure on $X$
which is absolutely continuous with respect to $\mu$. We assume moerover
that
\begin{itemize}
\item(bounded support)  there exists $r_1>0$ such that $\nu_x$ is supported in $B(x,r_1)$ for all $x\in X$, 
 
\item(non-degeneracy) there exists $c>0$ such  
that $p_x(y):=d\nu_x/d\mu(y)\geq c$  for all $y\in SB(x,r_0)$.
\end{itemize}

Denote by $P$ the Markov operator on $L^2(X)$ defined by
$$Pf(x)=\int f(gy)d\nu_x(y).$$ 

\begin{defn}
Under the previous assumptions, we call $(X,P)$ a geometric $G$-random walk. Moreover, if $P$ is self-adjoint, then $(X,P)$ is
called a symmetric geometric $G$-random walk.
\end{defn}
Denote
$dP_x(y)=d\nu_x(y)=p_x(y)d\mu(y)$ and $dP^n_x(y)=p^n_x(y)d\mu(y)$.

\begin{thm}\label{localfieldthm2}
Let $G$ be a unimodular elementary solvable group with exponential
growth. Then for every symmetric quasi-$G$-transitive random walk
$(X,P)$, we have
$$\sup_{x\in X}p_x^n(x)\approx e^{-n^{1/3}}.$$
\end{thm}
\bpr Since $X$ and $G$ are large-scale equivalent,
 Theorem~\ref{largescalethm} implies that their large-scale isoperimetric
profiles are asymptotically equivalent. Therefore $j_{X,2}\approx
\log t$. Theorem~\ref{random/sobolevThm} then implies that the
probability of return of any reversible random walk at scale large
enough, decreases like $e^{-n^{1/3}}$. 
To apply this to our random
walk $P$, we just need to check that for $k$ large enough, $P^k$
satisfies the conditions (i) to (iii) of Section \ref{probaSection}, with $h$ as big as we want. This follows easily from the definition of $r_0$, $r_1$ and $c$, and the fact that $S$ generates $G$. \epr

\bigskip
\footnotesize

\end{document}